\documentclass[conference]{IEEEtran}
\IEEEoverridecommandlockouts
\usepackage{cite}
\usepackage{amsmath,amssymb,amsfonts}
\usepackage{algorithmic}
\usepackage{graphicx}
\usepackage{subcaption}

\usepackage[english]{babel}
\usepackage{amsthm}

\theoremstyle{definition}

\theoremstyle{remark}

\usepackage{textcomp}
\usepackage{xcolor}
\def\BibTeX{{\rm B\kern-.05em{\sc i\kern-.025em b}\kern-.08em
    T\kern-.1667em\lower.7ex\hbox{E}\kern-.125emX}}
\begin{document}

\title{Improved microgrid resiliency through distributionally robust optimization under a policy-mode framework\\
\thanks{The authors are with the Electricity Infrastructure and Buildings Division at PNNL, Richland, WA 99354, USA.}
}

\author{\IEEEauthorblockN{ Nawaf Nazir}
\IEEEauthorblockA{\textit{Energy and Environment Directorate} \\
\textit{Pacific Northwest National Laboratory}\\
Richland, USA \\
nawaf.nazir@pnnl.gov}
\and
\IEEEauthorblockN{Thiagarajan Ramachandaran}
\IEEEauthorblockA{\textit{Energy and Environment Directorate} \\
\textit{Pacific Northwest National Laboratory}\\
Richland, USA \\
thiagarajan.ramachandaran@pnnl.gov}
\and
\IEEEauthorblockN{Soumya Kundu}
\IEEEauthorblockA{\textit{Energy and Environment Directorate} \\
\textit{Pacific Northwest National Laboratory}\\
Richland, USA \\
soumya.kundu@pnnl.gov}
\and
\IEEEauthorblockN{Veronica Adetola}
\IEEEauthorblockA{\textit{Energy and Environment Directorate} \\
\textit{Pacific Northwest National Laboratory}\\
Richland, USA \\
veronica.adetola@pnnl.gov}
}

\maketitle

\begin{abstract}
Critical  energy  infrastructure  are constantly under  stress  due  to  the  ever  increasing  disruptions  caused  by wildfires, hurricanes, other weather related extreme events and cyber-attacks.  Hence  it  becomes  important  to  make  critical infrastructure   resilient   to   threats   from   such cyber-physical events. Such events are however hard to predict and numerous in nature and type, making it infeasible to become resilient  to  all  possible cyber-physical  event as such  an approach  would  make  the  system  operation  overly  conservative. Furthermore, distributions of such events  are  hard  to  predict  and  historical data  available  on such  events  is sparse.  To  deal  with  these  issues,  we present a policy-mode framework that enumerates and predicts the  probability  of  various  cyber-physical  events on top of a distributionally robust optimization (DRO) that is robust to the  sparsity  of the available historical data. The proposed  algorithm  is  illustrated on  an islanded microgrid example: a  modified IEEE 123-node  feeder with  distributed energy  resources (DERs) and energy storage.
\end{abstract}

\vspace{-2.5mm}
\section{Introduction}
Critical energy infrastructure has undergone significant changes in the past decades with the increased penetration of distributed energy resources (DERs) and other inverter-interfaced generation. While this has helped to achieve green-house gas emission reduction goals~\cite{di2018decarbonization}, it has also made the energy grid more vulnerable to breakdowns due to the uncertainty and variability in renewable energy generation~\cite{zio2013vulnerability,nazir2019convex}. Furthermore, the risk of weather related outages such as wildfires, hurricanes and other natural disasters have increased in the past years~\cite{NERC2021} and so have the risk of possible cyber-attacks~\cite{HR360}. 
This work provides a framework for improved resiliency of critical infrastructure systems, with an emphasis on microgrid resiliency. Microgrids are a group of flexible energy resources operating together locally as a single controllable entity~\cite{farrokhabadi2019microgrid}.
In our framework, we assess the risk associated with various cyber-physical events and based on the risk assessment operate the system under a particular policy-mode. This allows us to be resilient against the risks without being overly conservative. Then in order to deal with the sparsity of data available on such cyber-physical events, we develop a distributionally robust optimization (DRO) formulation that is robust to a range of disturbance distributions.

\textbf{Literature review/related work:}
Chance constraints constitute a means to provide certain guarantees on constraint satisfaction under uncertainty~\cite{sahinidis2004optimization} and have found applications in several domains including in power system optimization under uncertainty.
However, a major drawback of chance constraint formulations is that they only consider the probability of constraint violation and not the impact. In many critical infrastructure systems, minimizing the impact of uncertainties is far more significant. Conditional value at risk (CVaR) approaches do account for the risk in constraint violation~\cite{summers2015stochastic}, however, many of these methods require assumptions on the probability distribution of the uncertainty and convex reformulations exist only for a very small set of such distributions (e.g., Gaussian), which may not hold in practice.

Distributionally robust chance constraint problems have been well studied in the literature~\cite{calafiore2006distributionally}.
Even though, these methods consider the risk in constraint violations, they can often be overly conservative and furthermore still require assumptions on the underlying distribution (e.g., moments) to be accurate. In case the uncertainty distribution is not known beforehand, sample based approaches have been used in literature~\cite{cheng2019partial}. However, in case of rare events sampling based approaches cannot be applied directly, since they would require an impractical number of samples to yield reasonable solutions~\cite{barrera2016chance}, something which may not be available in practice. To deal with this issue, Wasserstein ambiguity sets have been considered in the literature~\cite{cherukuri2020consistency} where the chance and risk constraints are required to hold for a family of distributions within a distance (Wasserstein distance) from the observed realizations of uncertainty.

Even though these methods are distributionally robust, they can be overly conservative in certain situations, especially when there are a large possible number of events and sources of uncertainty. To deal with the large number of physical events that can disrupt critical infrastructure systems in a manner that is not overly conservative, we formulate a policy-mode based DRO formulation. Here we decide the mode in which the system is operated under from amongst a finite number of modes. The operational policy-mode in turn is decided based on the risk-prediction of various events. Several works in literature have developed risk-based approaches that predict the risk of outages~\cite{kezunovic2022data} and various events such as proactively de-energizing grid components in the context of extreme wildfire risk~\cite{rhodes2020balancing}. Such risk-prediction can often be based upon publicly available data, such as the wildfire risk data from the United States Geological Survey, in case of predicting events such as Public Safety Power Shutoff (PSPS)~\cite{kody2022optimizing}.\\
The main advantages of our approach is that it is robust to discrepancies in sampled distribution from historical data by utilizing the Wasserstein ambiguity set, while at the same time avoids being overly conservative by considering the risk-prediction of various events and based on the risk-prediction metrics, operates the system under a particular policy-mode considering only a subset of possible events.

The rest of the paper is organized as follows: Sec.\,II provides the modeling of the microgrid setup and it's various chance constraint and risk constraint formulations. Sec.\,III introduces the Policy-mode framework that considers the risk-prediction of various events. Sec.\,IV presents the distributionally robust formulation that utilizes the Wasserstein ambiguity set framework. Sec.\,V provides the simulation results that showcase the efficacy of the proposed approach and, finally, Sec.\,VI provides the conclusions.
\section{Modeling and Problem Formulation}
\subsection{Baseline optimization}
We consider a microgrid with various DERs, including solar PVs, diesel generators (DGs), and energy storage (ES) units. The distribution network power-flow can be expressed by the following linear equations in Branch Flow Model (BFM) format~\cite{gan2014convex}:
\begin{subequations}
\label{eq:pf_linear}
\begin{align}
    0&=W_{n,k}-W_{m,k}+(S_{n,m,k}Z_{n,m}^*+Z_{n,m}S_{n,m,k})\label{eq:P2_volt_rel}\\
0&=\text{diag}(S_{n,m,k})-\sum_{o:m\rightarrow o}\text{diag}(S_{m,o,k})+S^{\text{net}}_{n,k}\label{eq:P2_power_balance}\\
 P^{\text{net}}_{n,k}&=P^\text{pv}_{n,k}+P^\text{dg}_{n,k}+P^\text{es}_{n,k}-P^\text{load}_{n,k} \label{eq:P1_node_real_balance}\\
Q^\text{net}_{n,k}&=Q^\text{pv}_{n,k}+Q^\text{dg}_{n,k}+Q^\text{es}_{n,k}-Q^\text{load}_{n,k}\label{eq:P1_node_reactive_balance}\
\end{align}
\end{subequations}
where: \begin{align*}
    W_{n,k}:=V_{n,k}\cdot V_{n,k}^*\,,\quad I_{n,m,k}:=i_{n,m,k}\cdot i_{n,m,k}^*\,.
\end{align*}
and where $V_{n,k}$ denotes the multi-phase voltage vector at bus $n\in \mathcal{N}$, $i_{n,m,k}$ denotes the multi-phase current flowing from bus $n$ to bus $m$, $S_{n,m,k}$, $P_{n,m,k}$ and $Q_{n,m,k}$ denote the apparent, active and reactive branch power flow from $n\rightarrow m$, $Z_{m,n}$ be the impedance matrix of branch $(n,m)\in \mathcal{L}$, $S_{n,k}^{\text{net}}=P_{n,k}^{\text{net}}+jQ_{n,k}^{\text{net}}$ be the net apparent power injection at bus $n$ and where $\mathcal{N}$ denotes the set of buses, $\mathcal{L}\subset \mathcal{N}\times \mathcal{N}$ denotes the set of all branches and $k$ denotes the time-step. 
In addition the following operational reliability and safety constraints are imposed:
\begin{subequations}
\label{eq:pf_constraints}
\begin{align}
    &\left|\text{diag}(S_{n,m,k})\right|\leq \overline{S}_{n,m} \label{eq:P1_line_const}\\
&\underline{V}^2\leq \text{diag}(W_{n,k})\leq \overline{V}^2 \label{eq:P1_volt_const}\\
& \left|S^\text{pv}_{n,k}\right|\leq \overline{S}^\text{pv}_n\,,\quad \left|S^\text{dg}_{n,k}\right|\leq \overline{S}^\text{dg}_n\,, \label{eq:gen_constraints}
\end{align}
\end{subequations}
The full description of these constraints and other constraints on battery State-of-charge (SoC) and solar and load curtailment is omitted here for lack of space and can be found in~\cite{nazir2022optimization}. 
Based on the above equations, the baseline optimization problem for the microgrid can be formulated as:
\begin{subequations}\label{eq:Base_Opt}
\begin{align}
    \min_x f(x)\\
    \text{s.t.} \ \ g(x)\le 0
\end{align}
\end{subequations}
where $x$ denotes the decision variables (e.g., DG dispatch), $f(x)$ is the problem objective (e.g., minimize operating cost) and $g(x)$ are the system constraint (e.g., line limits). An example of the objective function that can be employed here is given by:
\begin{align}
    f(x):=\sum_k \sum_n \left(c_1 \mathbf{P}^\text{dg}_k+c_2 \mathbf{P}^\text{es}_k+c_3 \mathbf{P}^\text{load}_k\right)
\end{align}
The optimization problem in~\eqref{eq:Base_Opt} does not consider uncertainty in system parameters or forecasts, which are critical for practical applications. The next section will consider the optimization problem under uncertainty. 

\subsection{Optimization under uncertainty}
Uncertainty in the system can arise from various adversarial cyber-physical events, like a DG outage, a load-masking attack or a sudden loss in solar generation. In order to be resilient against such events, it is important to have sufficient reserves in the dispatchable resources (DG, PV, storage, and interruptible load) so that these can adjust their output in response to adversarial events. The purpose of this work is to formulate probabilistic methods that determine the amount of reserves requires to deal with such adversarial events in a resilient and cost-effective manner. If we consider a bounded set of uncertainty $w\in \Omega_{\text{W}}$, then the stochastic optimization problem can be expressed as:
\begin{subequations}\label{eq:Stoch_Opt}
\begin{align}
    &\min_x f(x,w)\\
    \text{s.t.} \ \ &g(x,w)\le 0, \ \forall w \in \Omega_{\text{W}}
\end{align}
\end{subequations}
where $x$ also includes additional decision variables in the form of reserves.
However, the above optimization problem is infinite dimensional and cannot be solved in a scalable manner. In the next section, we briefly describe some of the traditional probabilistic methods that are used to solve~\eqref{eq:Stoch_Opt}.
\subsection{Chance constraint formulations}
In this section, we briefly present various stochastic optimization reformulations to solve~\eqref{eq:Stoch_Opt}. One common approach is to robustify~\eqref{eq:Stoch_Opt} against a specified set of adversarial events. In the most general form, the robust optimization problem can be expressed as:
\begin{subequations}\label{eq:arbit_opt}
   \begin{align}
   & \min_{x} \max_{{w}} f({x},w)\\
    s.t. & \ g(x,{w}) \le 0 \qquad \forall {w}\in \Omega_W
\end{align}
\end{subequations}
The above optimization problem can be reformulated using the explicit maximization method, details of which can be found in~\cite{bai2015robust}. However, a robust approach can be overly conservative, especially when the set $w\in \Omega_{\text{W}}$ is large. In such situations, chance constraint based methods can be employed that avoid constraint violation with a certain probability. A general chance constraint optimization can be expressed as:
\begin{subequations}\label{eq:CC_Opt}
\begin{align}
   & \min_x \mathbb{E}[f(x,w)]\\
    \text{s.t.} \ \ & \mathbb{P}[g(x,w)\le 0]\ge 1-\rho
\end{align}
\end{subequations}
where $\mathbb{E}$ is the expected value operator, $\mathbb{P}$ is the probability operator, and $\rho$ is the pre-specified allowable constraint violation. 
One of the main drawbacks of chance constraint methods is that it considers the probability of an event and not the magnitude of constraint violation (risk) of the event. In order to consider the risk of constraint violation, the Conditional value at risk or CVaR method has been employed in the literature~\cite{summers2015stochastic}, especially to deal with low probability high impact events.
One such CVaR formulation is obtained by utilizing the Markov generating function which can be expressed as:
\begin{subequations}\label{eq:CVaR_Opt}
\begin{align}
   & \min_x \mathbb{E}[f(x,w)]\\
    \text{s.t.} \ \ & \mathbb{E}[g(x,w) + t]_+\le t\rho
\end{align}
\end{subequations}
where $t$ is an optimization variable and $[.]_+=\max (.,0)$. 
Even though CVaR related methods account for low probability high impact events, they have an underlying assumption on the probability distribution of the uncertainty, i.e., the historical data used to represent $w \in \Omega_{\text{W}}$ matches the actual distribution. However, in many practical situations the historical data on adversarial events can be very sparse and hence may not represent the true distribution of such adversarial events. 
In such cases, distributionally robust optimization (DRO) methods have shown promise~\cite{cherukuri2020consistency}. DRO methods such as the ones based on ambiguity balls consider a family of probability distributions within a certain metric in the optimization. By changing the radius of the ambiguity balls, the DRO formulation can be made more robust to lack of sufficient data to represent the disturbance distribution. However, these DRO methods do not deal with the potential conservativeness of such DRO when the adversarial set is large. The contribution of this paper is a Policy-mode based framework that reduces the conservativeness of DRO while at the same time accounting for the disturbances in a distributionally robust manner. The Policy-mode based framework achieves this by selecting a subset of disturbance events from a larger pool based on predicted risks. In the next section, we present the Policy-mode based framework in detail.
\vspace{-2.5mm}
\section{Policy Mode Framework}
\begin{figure}[t]
\centering
\includegraphics[width=0.45\textwidth]{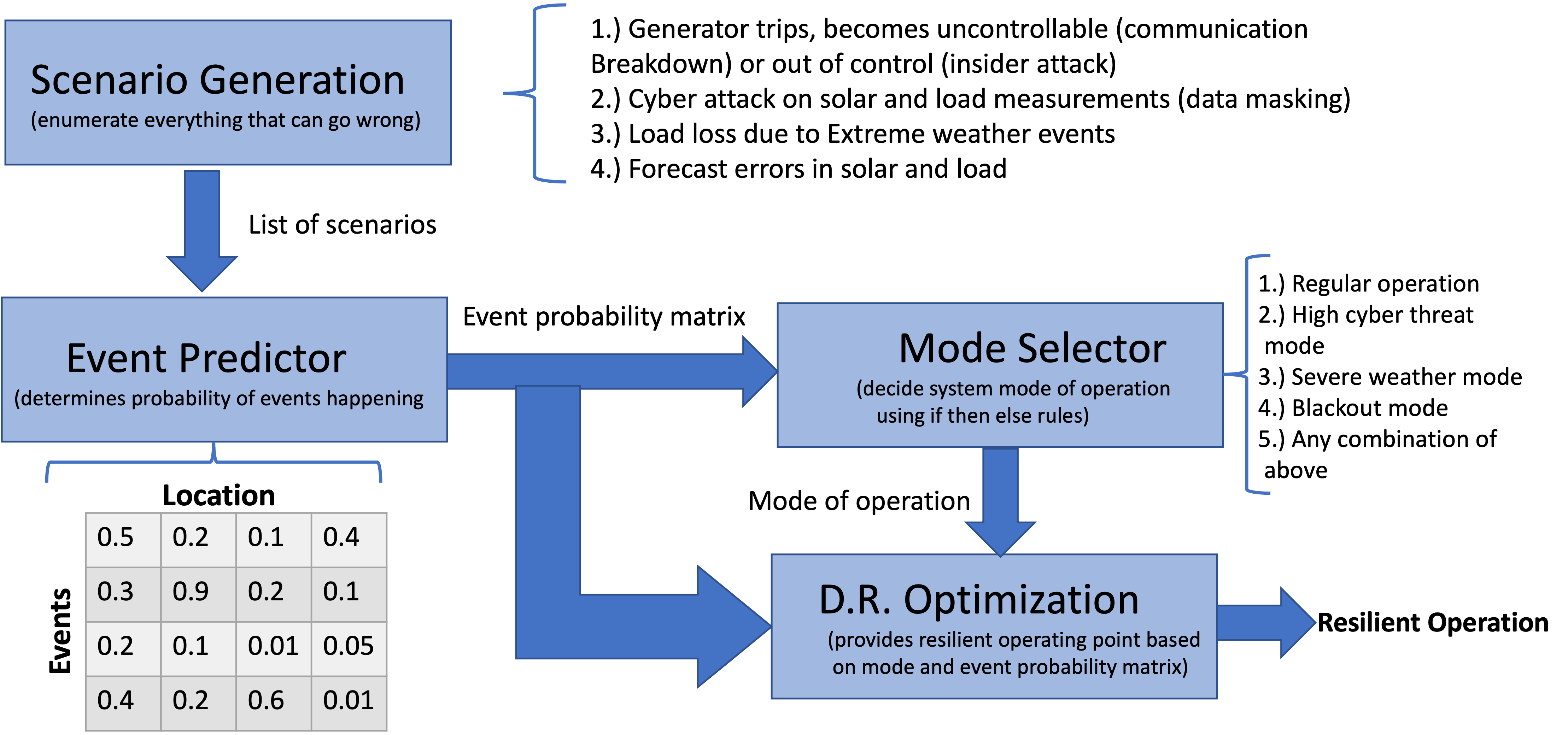}
\caption{\label{fig:Event_pred}Policy mode implementation framework.}
\end{figure}

\vspace{-2.5mm}
In this section we present a Policy-mode framework implementation for the DRO formulation which is also depicted in Fig.~\ref{fig:Event_pred}.
The entire framework comprises of four building blocks, each of which will be described in the ensuing text.

\textbf{Scenario Generation:}
In the first stage of the framework,
we enumerate the different events, adversarial attacks, breakdowns, faults and forecast error that can happen in the system and their location. These can be collected based on the history of past events.

\textbf{Event Prediction:}
After the scenario generation, the list of possible events collected is passed on to an Event Predictor block. The role of the event predictor is to determine the probability of various events occurring at different locations.
Examples of methods that could be utilized for such event prediction exist in literature for outages caused by wildfires and other weather related events~\cite{kezunovic2022data,rhodes2020balancing,kody2022optimizing}.
The output of the Event Predictor block is an Event-Location probability matrix that predicts the probability of an event occurring at a particular location as illustrated in Fig.~\ref{fig:Event_pred}.

\textbf{Mode Selection:}
After the Event prediction block, the obtained Event Probability matrix is fed into the Mode Selector block. The role of the Mode Selector is to reduce the dimensionality of the problem, as the number of possible events to consider can be very large, making the optimization problem intractable.
As an example, Fig.~\ref{fig:Event_pred} lists the modes that the system can operate in and depending upon the values in the Event probability matrix and their impact, the Mode Selector block chooses to operate the system in one or a combination of these modes.

\textbf{D.R. Optimization:}
Finally, the modes of operation along with the Event probability matrix is fed to the D.R.Optimization block, that solves a DRO problem and provides resilient operating points and reserves required by the system to operate. Details about the DRO formulation will be provided in the next section.


\section{Distributionally Robust Optimization Formulation}
In this Section we summarize the Distributionally robust optimization (DRO) formulation based on the Wasserstein ambiguity set that has previous been presented in several works in literature including~\cite{duan2018distributionally,zhou2020linear,shi2018distributionally}.

\subsection{Wasserstein Ambiguity Set}

The Wasserstein distance $W(.,.)$ between the true distribution $\mathbf{P}$ and sample probability distribution $\mathbf{\hat{P}}_N$ is used to define the ambiguity set $\hat{P}_N$ given by:
\begin{align}\label{eq:Wass-exp}
    \hat{P}_N=\{P\in P(\Xi): W(\mathbf{P,\mathbf{\hat{P}}_N}) < \epsilon(N)\}
\end{align}
where $P(\Xi)$ denotes the set of all probability distributions with support $\Xi$ and where:
\begin{align}
    \mathbf{\hat{P}}_N=\frac{1}{N}\sum_{k=1}^N \delta _{\hat{\zeta}^{(k)}}
\end{align}
where $\hat{\zeta}^{(k)} \in \mathbb{R}^N$ are the samples of distribution and $N$ is the number of such samples available.
Based on the this, the stochastic constraint of the form in~\eqref{eq:CC_Opt}, can be expressed in the DRO form as:
\begin{align}\label{eq:Wass-stoch}
\inf_{\mathbb{P}\in \hat{P}_N} \mathbb{P}[g(x,w)\le 0]\ge 1-\rho
\end{align}
Since the above constraint is intractable, the goal is to develop a convex approximation of~\eqref{eq:Wass-stoch} of the form:
\begin{align}\label{eq:set_U}
    g(x,w)\le 0, \ \forall w \in \mathcal{U}
\end{align}
where $\mathcal{U}$ is a deterministic uncertainty set such that the robust deterministic constraint in~\eqref{eq:set_U} implies satisfaction of the DRO constraint in~\eqref{eq:Wass-stoch}. 
Based on the formulation in~\cite{duan2018distributionally}, the set $\mathcal{U}$ can be expressed as:
\begin{align}
    \mathcal{U}=Conv(\{u^{(1)}, u^{(2)}, \hdots u^{(2m)}\})\\
    u^{(i)}=\hat{\Sigma}^{1/2}v^{(i)}+\hat{\mu}, \ 1\leq i \leq 2^m
\end{align}
where $Conv$ represents the convex combination of the terms, $v=\hat{\Sigma}^{-1/2}(w-\hat{\mu})$, where $\hat{\mu}$ is the sample mean and $\hat{\Sigma}$ is the sample co-variance.
This results in the following deterministic DRO formulation based on Wasserstein ambiguity set:
\begin{subequations}\label{eq:DRW_Opt}
\begin{align}
&\min_x f(x)\\
\text{s.t.} \ \ &g(x,u^{(i)})\leq 0, \ 1\leq i \leq 2^m
\end{align}
\end{subequations}
Since the set $\mathcal{U}$ can be calculated beforehand, the optimization problem~\eqref{eq:DRW_Opt} can be solved in a tractable manner using off the shelf convex optimization solvers. 

It should be noted that the robustness of the DRO formulation can be tuned through the parameter $\epsilon$, in a similar manner to the allowable constraint violation parameter $\rho$ in case of CVaR formulation. However, the difference between the two is that DRO accounts for the robustness within the distribution of the disturbance events itself by considering a family of distributions through some metric, whereas CVaR assumes a particular distribution for the disturbance events and only provides performance guarantees within that assumed distribution.
However, DRO can be overly conservative, especially when the possible number of disturbance events are large. This is where the policy-mode framework is able to reduce the dimensionality of considered events by operating the system in one of a certain finite number of operating modes.

The next section will provide simulation results on an IEEE-123 node microgrid system that showcases the efficacy of the developed DRO formulation and also compares it with the CVaR formulation.

\section{Simulation Results}
To illustrate the DRO formulation, we peforms simulations on a modified three-phase IEEE-123 node systems~\cite{kersting2006distribution} converted into an islanded microgrid by disconnecting the substation from the utility. Further details of the simulation setup are omitted due to space constraints and can be found in~\cite{nazir2022optimization}. Table~\ref{tab:scenarios} lists all the possible disturbance events along with their location in the system.
\begin{table}[htpb]
    \centering
    \begin{tabular}{c|c|c}
         \textsc{\#} & \textsc{Events} &  \textsc{Disturbance Location}\\
         \hline\hline
         1 & DG trip & \{47, 49, 65, 76\}\\
         \hline
         2 & PV trip & \{28, 42, 60, 76, 100\}\\
         \hline
         3 & DG cyber-attack & \{47, 49, 65, 76\}\\
         \hline
         4 & PV cyber-attack & \{28, 42, 60, 76, 100\}\\
         \hline
         5 & Load cyber-attack & \{1, 6, 28, 29, 35, 42, 47, 48,\\ 
         && 49, 52, 60, 64, 65, 76, 77, 100\}\\
         \hline
         6 & PV forecast error & \{28, 42, 60, 76, 100\}\\
         \hline
         7 & Load forecast error & \{1, 6, 28, 29, 35, 42, 47, 48,\\ 
         && 49, 52, 60, 64, 65, 76, 77, 100\}\\
         \hline
         8 & Severe weather PV loss & \{28, 42, 60, 76, 100\}\\
         \hline
         9 & Severe weather Load loss & \{1, 6, 28, 29, 35, 42, 47, 48,\\ 
         && 49, 52, 60, 64, 65, 76, 77, 100\}
    \end{tabular}
    \caption{List of events and their locations}
    \label{tab:scenarios}
\end{table}
These list of scenarios are fed to the Event Predictor block that predicts the probability of various events. Next the Event probability matrix is fed to the Mode Predictor, which based on the impact of different events, decides to operate the system in a particular mode from amongst a finite number of modes. 
We solve the DRO with the event prediction risk probabilities based on the selected mode of operation. The output of the DRO provides the operating point of the generation and curtailment, along with the reserve values that would be required when the uncertainty reveals itself. 
To validate the resiliency of the optimization, we perform a series of simulations in GridLab-D~\cite{chassin2008gridlab} by sampling from the true uncertainty distribution.

As in our case the sampled distribution does not match the actual distribution, as depicted in Fig.~2(a), so the assumptions of the CVaR formulation break-down and as a result we observe large constraint violations as depicted in Fig.~2(b).
Next we apply the DRO formulation based on Wasserstein ambiguity set on the same distribution. As seen by the results in Fig.~2(d), the calculated reserves are able to manage the uncertainty in this case, as it accounts for the disparity between the sampled distribution and actual distribution.
Finally, to highlight the trade-off between resiliency (available reserves) and efficiency (inverse of cost), Fig.~2(c) shows the trade-off under 4 different modes of operation, when using  DRO, highlighting the difference across different operation modes.

\begin{figure}
    \setkeys{Gin}{width=\linewidth}
    \begin{subfigure}[t]{0.24\textwidth}{(a)}
           \includegraphics{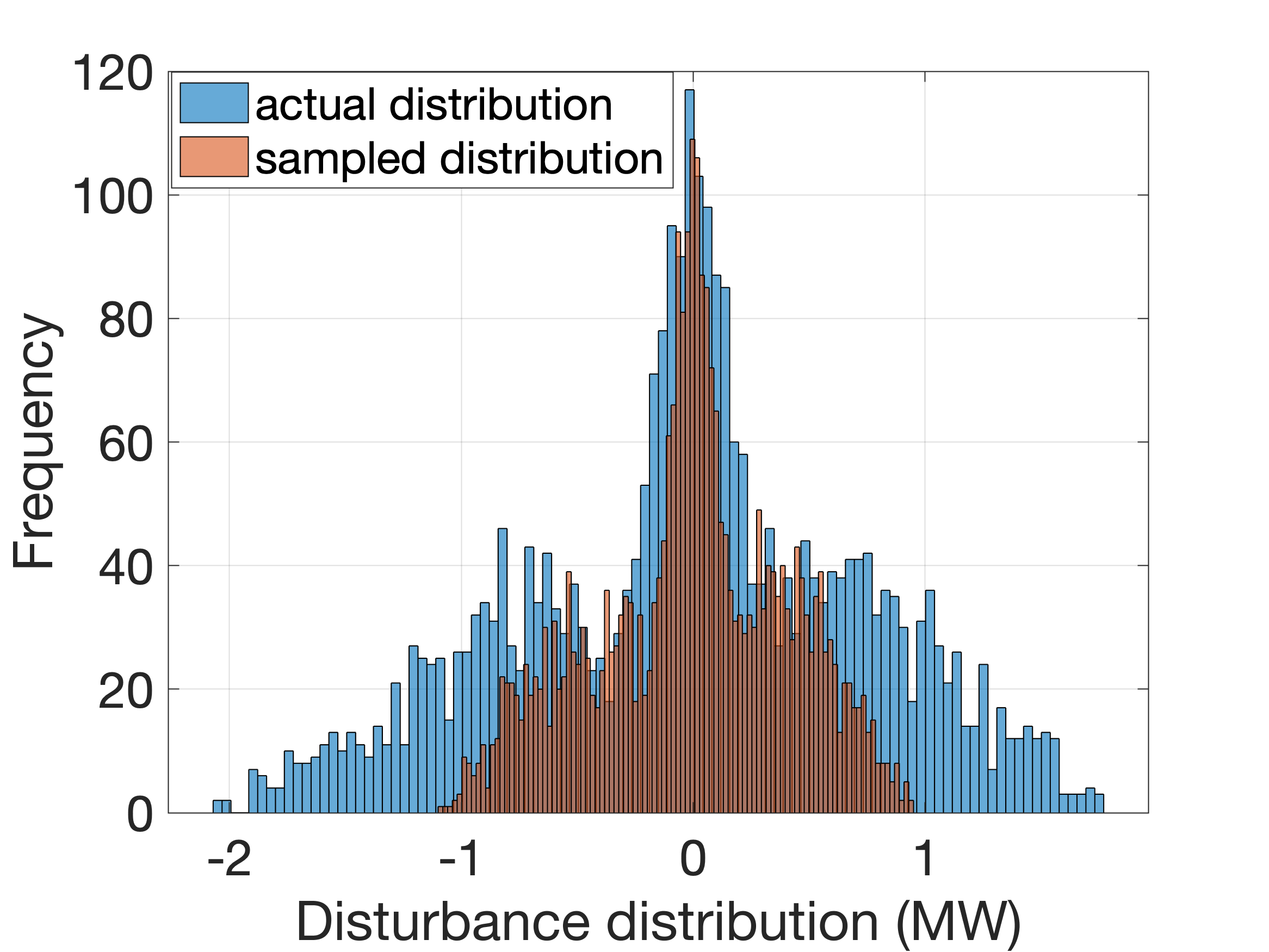}
        \label{fig:unmatched_dist_comp}
    \end{subfigure}
    \hfill
    \begin{subfigure}[t]{0.24\textwidth}{(b)}
        \includegraphics{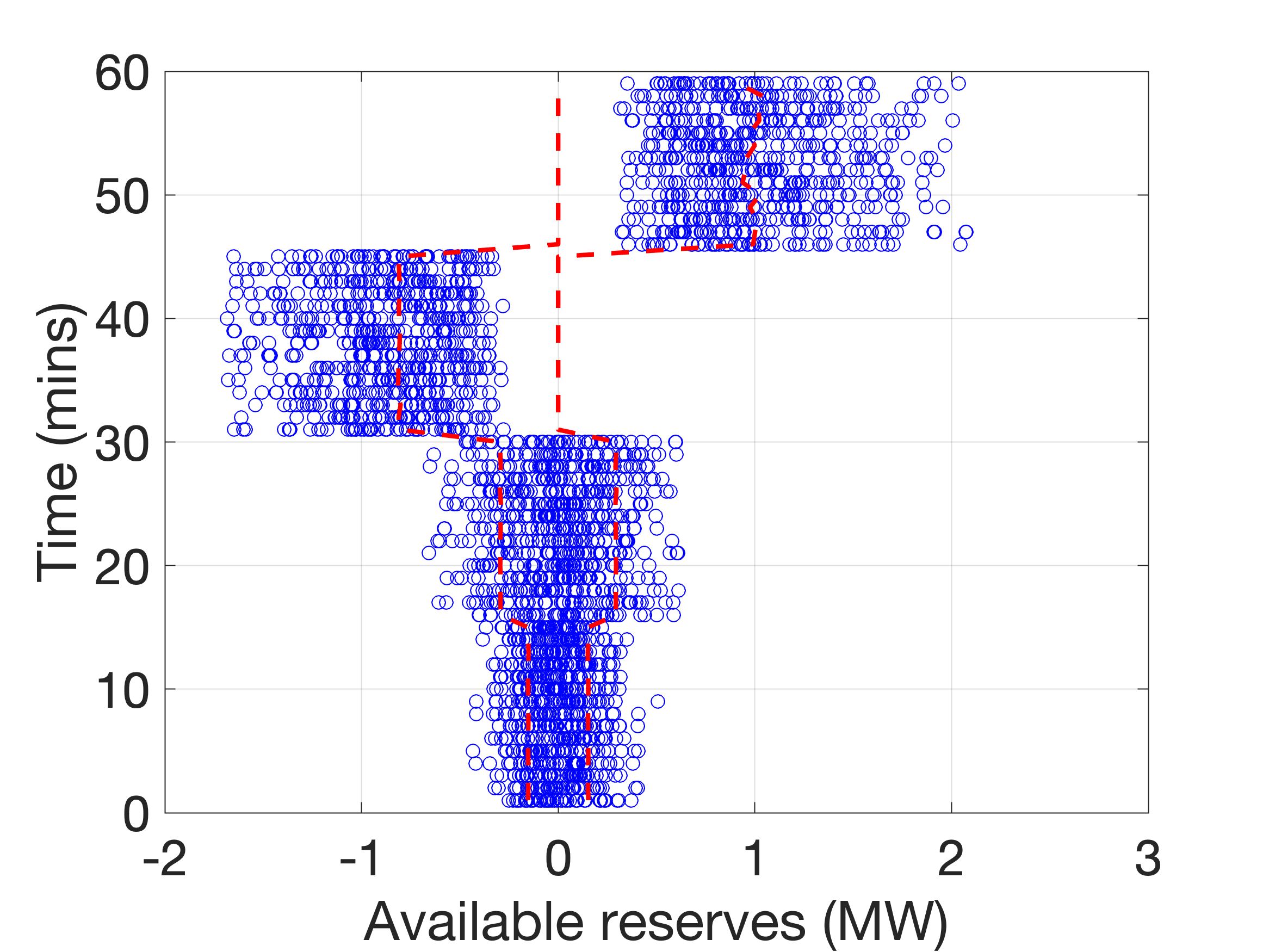}
\label{fig:CVaR_unmatched}
    \end{subfigure}
    \\
    \begin{subfigure}[t]{0.24\textwidth}{(c)}
    \includegraphics{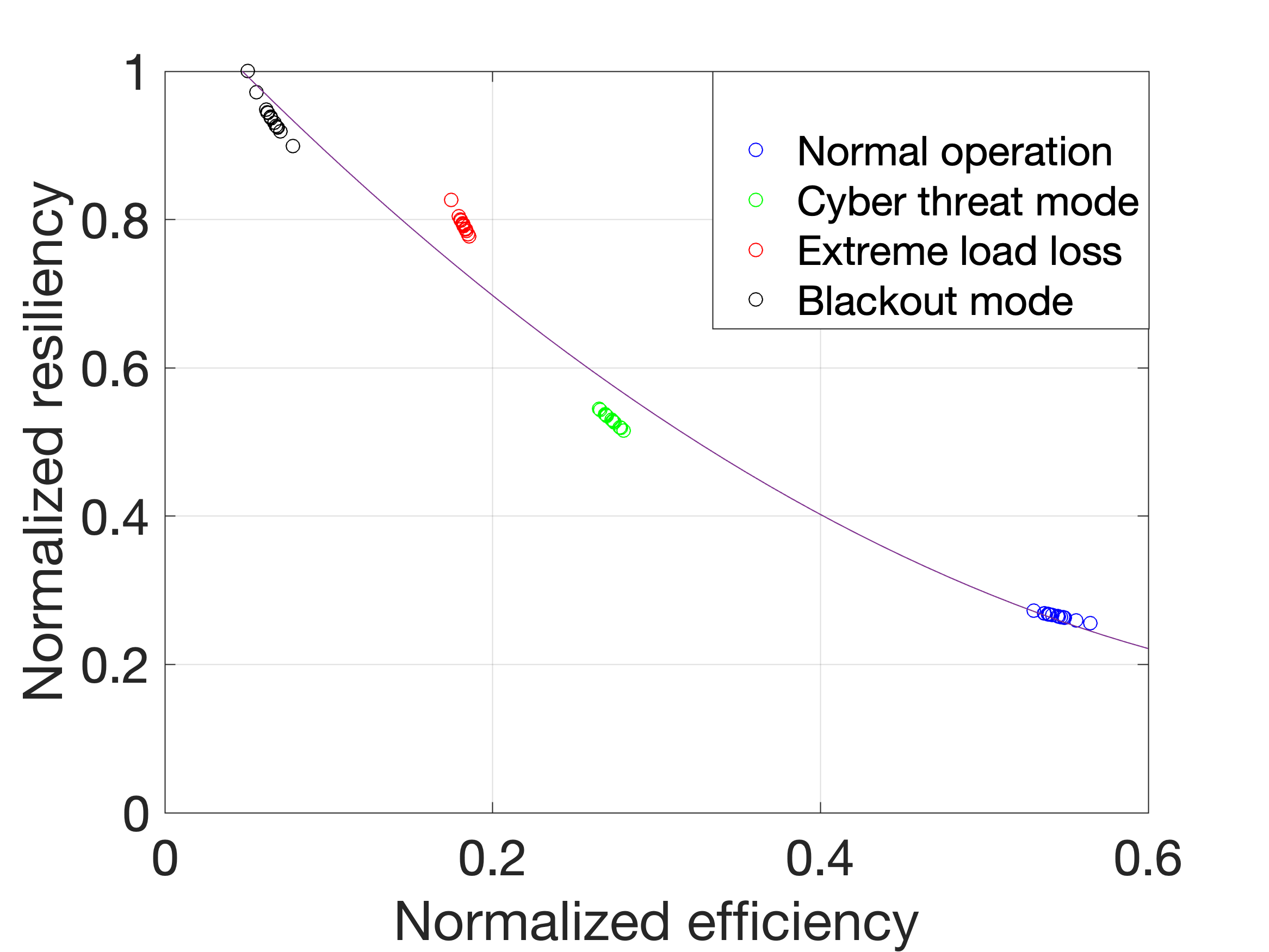}
    \label{fig:eff_res_DRO}
\end{subfigure}
\hfill
\begin{subfigure}[t]{0.24\textwidth}{(d)}
    \includegraphics{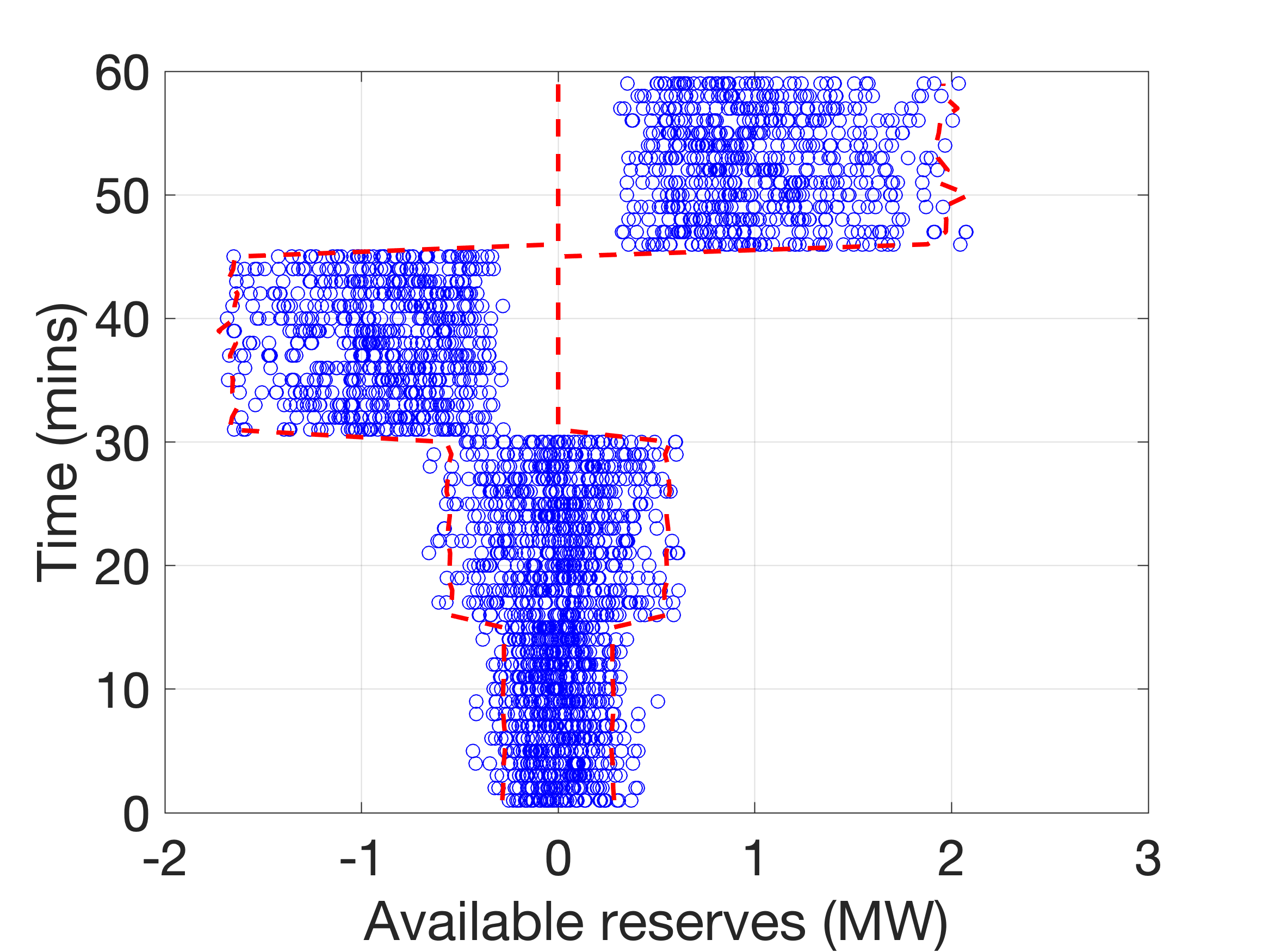}
\label{fig:DRW_unmatched}
\end{subfigure}
    \caption{(a) Comparison of histograms of sampled and actual distribution (b) Disturbance impacts (blue dots) in comparison to available reserves (red lines) in CVaR case (c) Trade-off between efficiency and resiliency for DRO (d) Disturbance impacts (blue dots) in comparison to available reserves (red lines) in the DRO case.}
\end{figure}
\section{Conclusion}
\vspace{-1mm}
This paper presented a Policy-mode framework based on available risk-predictions of various events to formulates a data-driven distributionally robust optimization (DRO) for a microgrid. Wasserstein ambiguity sets are utilized to deal with the sparsity in available historical data. Through simulation results we showcase the efficacy of our approach as compared to other risk based methods. 


\section*{Acknowledgment}
This research was supported by the ``Resilience through Data-driven Intelligently-Designed Control'' (RD2C) Initiative, under the Laboratory Directed Research and Development (LDRD) Program at PNNL. PNNL is a multi-program national laboratory operated for the U.S. Department of Energy by Battelle Memorial Institute under Contract No. DE-AC05-76RL01830. 


\bibliographystyle{IEEEtran}
\bibliography{References.bib}

\end{document}